\documentclass[12pt]{article}
\usepackage[top=1in,bottom=1in,left=1in,right=1in]{geometry}
\usepackage{indentfirst}
\usepackage{latexsym}
\usepackage{graphicx}
\usepackage{mathrsfs}
\usepackage{amsfonts,amsmath,amsthm,amssymb,mathtools}
\usepackage{latexsym, euscript, epic, eepic, color}

\newtheorem{lemma}{Lemma}[section]
\newtheorem{theorem}{Theorem}[section]

\begin{document}
\baselineskip=19pt

\newtheorem{conjecture}[theorem]{Conjecture}

\title{Domination number in block designs}
\author{Lang Tang, Shenglin Zhou\footnote{Corresponding author. This work is
supported by the National Natural Science Foundation
of China (Grant No.11471123). slzhou@scut.edu.cn}\\
\small \it School of Mathematics, South China University of Technology,\\
\small \it Guangzhou 510640, P.R. China\\
\date{}
}
\maketitle

\begin{abstract}
Let $G=(V,E)$ be a simple connected graph. A set of vertices $S\subseteq V$ is said to be a dominating set if for any vertex in $V\setminus S$ is
adjacent to at least one vertex in $S$. The domination number $\gamma(G)$ of $G$ is the minimum cardinality among all such sets.

In this paper, we obtain some results on the domination number of the incidence graphs of combinatorial designs. In particular, we prove a conjecture and disprove another conjecture in a recent paper by Goldberg, Rajendraprasad and Mathew. We also prove a third conjecture by the same authors for block-transitive symmetric designs.
\end{abstract}


\medskip
\noindent{\bf MR(2000) Subject Classification:} 05C69, 05B05, 05E18

\section{Introduction}

Let $G=(V,E)$ be a simple, undirected graph, where $V$ is the vertex set of
$G$ and $E$ is the edge set of $G$. Let $S\subseteq V$ be a set of vertices. A vertex $u\in V$ is said to be dominated by $S$ if either $u\in S$ or $u$ is adjacent to at least one vertex in $S$. $S$ is said to be a \textit{dominating set} if every vertex of $G$ is dominated by $S$. The \textit{domination number} $\gamma(G)$ of $G$ is the minimum cardinality among all dominating sets of $G$. A dominating set of $G$ whose cardinality is $\gamma (G)$ is called a  \textit{minimum dominating set} of $G$.

Hedetniemi and Laskar (1990) noted in \cite{Hedetniemi} that the problem of domination number can be dated back to at least the 1950's,
by K\"{o}nig, Berge and Ore et al., and has big advancement in the middle 1970's.
Recently, researches on all kinds of dominating sets in graphs and relationships between domination and other graphic parameters have become an very
important field in graph theory. For a general graph, the problem of finding its minimum domination set is a NP-hard problem \cite{Garey}.

Let $v$, $k$, $\lambda$ be positive integers such that $v\geq k\geq 2$. Let $X=\{x_1,x_2,\ldots,x_v\}$ be a finite set of $v$ elements, called points, and $\mathcal{B}$ be a family of $k\text{-}$subsets, called blocks. The pair $D=(X,\mathcal{B})$ is called a \textit{$2\text{-}(v,k,\lambda)$ block design}, or simply a \textit{$2\text{-}$design}, if every pair of distinct points are contained in exactly $\lambda$ blocks. The integers $v,k$ and $\lambda$ are called parameters of $D$. A pair $(x,B)$ with $x\in X$ and $B\in \mathcal{B}$ is called a \textit{flag} of $\mathcal{D}$ if $x\in B$. The set of all flags of $D$ is denoted by $\mathcal{F}$. Let $b=|\mathcal{B}|$ denote the total number of blocks. By \cite{Wan}, every point occurs in exactly $r=\frac{\lambda(v-1)}{k-1}$ blocks, $bk=vr$ and $b\geq v$. For $x\in X$, let $B(x)=\{B_1,B_2,\ldots,B_r\}$ be the family of blocks containing $x$, called the \textit{pencil} of $x$. A design $D=(X,\mathcal{B})$ is called \textit{a symmetric
design} if $b=v$.
 The
\textit{incidence graph} of $D$ is defined by $G_D=(V(G_D),E(G_D))$, where $V(G_D)=X\cup \mathcal{B}$ and $E(G_D)=\mathcal{F}$.

Both combinatorial designs and domination in graphs have been widely studied. While they have been hardly attempts to marry this two subjects. Laskar and
Wallis have obtained some results about the domination number of the line graph of $G_D$ in \cite{Laskar}.  However, they did not consider the domination
number of $G_D$ itself. Goldberg, Rajendraprasad and Mathew attempted to combine this two subjects first time in \cite{Felix} and got some interesting results. In this paper, our main purpose is to study the domination number of designs and solve the following three conjectures proposed in \cite{Felix}.
\begin{conjecture}\label{finite} {\rm\cite [Conjecture 6.2]{Felix}}
Finite projective planes are super-neat.
\end{conjecture}

\begin{conjecture}\label{symmetric}{\rm \cite[Conjecture 8.6]{Felix}}
Let $D=(X,\mathcal{B})$ be a symmetric $2$-$(v,k,2)$ design with $k\geq 4$. Then $\gamma(D)=k.$
\end{conjecture}

\begin{conjecture}\label{residual}{\rm \cite[Conjecture 9.2]{Felix}}
Let $D=(X,\mathcal{B})$ be a symmetric $2\text{-}(v,k,\lambda)$ design and $D_1=Res(X,\mathcal{B},B_0)$, a residual design of $D$, where $B_0\in \mathcal{B}$.
Then
$\gamma(D_1)=\gamma(D)-1.$
\end{conjecture}

In this article, we first study the domination number of $2\text{-}(v,k,\lambda)$ designs and non-symmetric $2\text{-}(v,k,1)$ designs. Our main results are as follows.
\begin{theorem}\label{g}
Let $D=(X,\mathcal{B})$ be a $2\text{-}(v,k,\lambda)$ design. Then
$\gamma (D)\geq \left\lceil\frac{2v-1-\frac{k-1}{\lambda}}{k} \right\rceil.$
\end{theorem}

\begin{theorem}\label{NO6}
Let $D=(X,\mathcal{B})$ be a non-symmetric $2\text{-}(v,k,1)$ design. Then
\[2k-1\leq \gamma (D)\leq \frac{(v-k^2)(v-1)}{k(k-1)}+2k-1.\]
\end{theorem}

We then give a sufficient condition for super-neat designs, which helps us prove Conjecture \ref{finite}. Besides, we obtain the domination number of affine planes and prove that affine planes are super-neat.

The next result shows that Conjecture \ref{symmetric} is wrong in general.

\begin{theorem}\label{NO8}
Let $D=(X,\mathcal{B})$ be a symmetric $2\text{-}(v,k,2)$ design with $k\geq 5$. Then
$$\gamma(D)\geq k-1+\sum_{i\geq 1}^{}\left\lfloor \frac{k-4}{2^{2i-1}} \right\rfloor.$$
\end{theorem}

\begin{theorem}\label{NORes}
Let $D=(X,\mathcal{B})$ be a block-transitive symmetric $2\text{-}(v,k,\lambda)$ design and $D_1=Res(X,\mathcal{B},B_0)$, a residual design of $D$, where $B_0\in
\mathcal{B}$. Then
$\gamma(D_1)=\gamma(D)-1.$
\end{theorem}
This theorem proves Conjecture \ref{residual} in the case when $D$ is a block-transitive symmetric $2\text{-}(v,k,\lambda)$ design.

\section{Preliminaries}

In this section, we introduce notations and give a few preliminary results which will be used throughout this paper. Undefined notations can be found in  \cite{Bondy}.

Let $G=(V,E)$ be a simple connected undirected graph. A \textit{bipartite graph} is one whose vertex set can be partitioned into two subsets $X$ and $Y$ such that every edge has one end in $X$ and the other one in $Y$; such a partition
$(X,Y)$ is called a bipartition of the graph. Let $S\subseteq V$ and $u\in S$. A vertex $v\in V\setminus S$ is called an \textit{external private neighbour} of $u$ if $u$ is the only neighbour of $v$ in $S$.

Let $D=(X,\mathcal{B})$ be a $2\text{-}(v,k,\lambda)$ design. It is easy to see that $G_D$ is a bipartite graph with bipartition $(X,\mathcal{B})$. 
From now on, we simply denote by $\gamma(D)$ the domination number $\gamma(G_D)$ of $G_D$. For a set of points $P$ of $D$, the blocks of $D$ are naturally partitioned into two parts, namely $L(P)=\{B\in \mathcal{B}|B\cap P\neq \emptyset\}$ and $\hat{L}(P)=\{B\in \mathcal{B}|B\cap P=\emptyset\}$. Denote $I_P=P\cup \hat{L}(P)$. For a subset $S\subseteq V(G_D)$ of $D$, we define $\pi(S)=S\cap X$. We say that $S$ is a \textit{neat set} if $S=I_P$ for some set of points $P$ of $D$. Obviously, $P=\pi(S)$. $D$ is said to be a \textit{neat design} if $D$ has a neat dominating set $S$ with $|S|=\gamma(D)$. If all minimum dominating sets of $D$ are neat, we say that $D$ is a \textit{super-neat design}.

\begin{lemma}\label{4.1}{\rm \cite{Felix}}
Let $D=(X,\mathcal{B})$ be a $2\text{-}(v,k,\lambda)$ design, $S$ be a dominating set of $G_D$ and $P=\pi(S)$. Then $\hat{L}(P)\subseteq S$.
\end{lemma}

By this Lemma, it is easy to see that $I_P\subseteq S$. This means that $S$ is not a neat set if and only if $I_P\subsetneq S$.

\begin{lemma}\label{4.4}{\rm \cite{Felix}}
Let $D=(X,\mathcal{B})$ be a $2\text{-}(v,k,\lambda)$ design, $S$ be a dominating set of $G_D$ and $P=\pi (S)$. Then
\[|S| \geq \left\lceil {\frac{{v+|P|(k-1)}}{k}}\right\rceil. \]
\end{lemma}

\begin{lemma}\label{5.2}{\rm \cite{Felix}}
Let $D$ be a finite projective plane of order $q$. Then $\gamma (D)=2q$.
\end{lemma}

\begin{lemma}\label{4.6} {\rm\cite{boll}}
Let $G$ be a graph without isolated vertices. Then $G$ has a minimum dominating set $S$ in which every vertex has an external private neighbour.
\end{lemma}


\section{Domination number in $2\text{-}$designs}

\textbf{Proof of Theorem \ref{g}.}
Let $S$ be a minimum dominating set of $G_D$ as provided by Lemma \ref{4.6} and $P = \pi(S)$. If $P=X=S$, then for any $B\in\mathcal{B}$, $S_1=\{B\}\cup (P\setminus B)$ is a dominating set of $G_D$ and $|S_1|=v-k+1<v=|S|$, a contradiction. Now we assume that there exists a point $x\in X\setminus P$. As $S$ is a dominating set, there exists a block $B_0\in S$ such that $x\in B_0$. By Lemma \ref{4.6}, $B_0$ has an external private neighbour $y \notin S$, that is $y \notin P$. Then all other $r-1$ blocks containing $y$ are not in $S$, which means that each of these blocks contains at least one point of $P$. On the other hand, every point $z$ in $P$ exactly dominates $\lambda$ blocks in $B(y)\setminus \{B_0\}$ since $\{z,y\}$ contained in exactly $\lambda$ blocks. Hence $|P|\geq \left\lceil\frac{r-1}{\lambda}\right\rceil \geq \frac{r-1}{\lambda}$.

By Lemma \ref{4.4} and $r(k-1)=\lambda (v-1)$, we have
\[\gamma(D)=|S|\geq \frac{v+|P|(k-1)}{k}\geq\frac{v+\frac{r-1}{\lambda}(k-1)}{k}=\frac{2v-1-\frac{k-1}{\lambda}}{k}.\]
This proves Theorem \ref{g}.\qed

\medskip
To prove Theorem \ref{NO6}, we need the following lemma.

\begin{lemma}\label{NO1}
Let $D=(X,\mathcal{B})$ be a $2\text{-}(v,k,\lambda)$ design and $P \subseteq X$ with $|P| \leq \left\lceil{\frac{r}{\lambda }} \right\rceil-1$ , then $I_P$ is a dominating set of $G_D$.
In particular, when $\lambda=1$, then $I_P$ is a dominating set of $G_D$ if $|P| \leq r-1$.
\end{lemma}
\noindent\textbf{Proof.} Assume that $I_P$ is not a dominating set. Then there exists $x \in X\backslash P$ such that $B(x)\cap \hat{L}(P)=\emptyset$, that is  $B(x)\subseteq L(P)$. Then for any $B_{i}\in B(x)$, there exists $x_{i}\in P\cap B_{i}$, $i=1,2,\ldots,r$,
and every point $y$ in $P$ exactly dominates $\lambda$ blocks in $B(x)$ for $\{x, y\}$ contained in exactly $\lambda$ blocks.
Thus, $|P|\geq \left\lceil{\frac{r}{\lambda }} \right\rceil$, a contradiction.\qed

\medskip
\noindent\textbf{Proof of Theorem \ref{NO6}.}
For any $B_0\in \mathcal{B}$, let $P=B_0$, then $|P|=k<r$. It follows from Lemma \ref{NO1} that $I_P$ is a dominating set of $G_D$. Since $|L(P)|=(r-1)k+1$, then
\begin{eqnarray*}
|\hat{L}(P)|=b-|L(P)|=b+k-kr-1=\frac{(v-k^2)(v-1)}{k(k-1)}+k-1.
\end{eqnarray*}
Hence
\[\gamma (D)\leq |S|=|P|+|\hat{L}(P)|=\frac{(v-k^2)(v-1)}{k(k-1)}+2k-1.\]

On the other hand, since $D$ is a non-symmetric $2\text{-}(v,k,1)$ design, we have $v\geq k^2$. By Theorem \ref{g}, we then have
\[\gamma(D)\geq \frac{2v-1-\frac{k-1}{\lambda}}{k}\geq 2k-1.\]
This completes the proof.\qed

\medskip
Now we proceed to prove Theorem \ref{NO8}.

\begin{lemma}\label{NO7}
Let $D=(X,\mathcal{B})$ be a symmetric $2\text{-}(v,k,2)$ design and $P\subseteq X$. If $2\leq |P|=\ell \leq k$, then
\[|L(P)|\leq \frac{\ell}{2}(2k-1-\ell)+1\]
and
\[|I_P|\geq \frac{k^2-k}{2}+\frac{\ell^2-\ell(2k-3)}{2}.\]
Moreover, the equalities hold if and only if there exists $B_0\in \mathcal{B}$ such that $P\subseteq B_0$.
\end{lemma}
\noindent\textbf{Proof.} We proceed by induction on $\ell$. If $\ell=2$, then $|L(P)|=2k-2$. Thus
\[|I_P|=\ell+v-|L(P)|=\frac{k^2-k}{2}+\frac{2^2-2(2k-3)}{2},\]
and there exists $B_0\in \mathcal{B}$ such that $P\subseteq B_0$ as $\ell=2$.

Assume that $|P|=\ell-1<k$ and the result holds. Then for any $x_0\in X\setminus P$, $|L(P)\cap L(\{x_0\})|\geq |P|+1=\ell$, and
the equality holds if and only if there exists $B_0\in \mathcal{B}$ such that $\{x_0\}\cup P \subseteq B_0$.

Let $P'=P \cup \{x_0\}$. Then $|P'|=\ell$ and
\begin{eqnarray*}
|L(P')|&=&|L(P'-\{x_0\})\cup L(\{x_0\})|\\
&=&|L(P)\cup L(\{x_0\})|\\
&=&|L(P)|+|L(\{x_0\})|-|L(P)\cap L(\{x_0\})|\\
&\leq &\frac{\ell-1}{2}(2k-1-\ell+1)+1+k-\ell\\
&=&\frac{\ell}{2}(2k-1-\ell)+1.
\end{eqnarray*}
The equality holds if and only if there exists $B_0\in \mathcal{B}$ such that $P'=\{x_0\}\cup P \subseteq B_0$.

Thus
$|L(P)|\leq\frac{\ell}{2}(2k-1-\ell)+1$
and
\[|I_P|=|P|+|\hat{L}(P)|\geq \frac{k^2-k}{2}+\frac{\ell^2-\ell(2k-3)}{2}\]
hold for all $2\leq |P|=\ell\leq k$. The equalities hold if and only if there exists $B_0\in \mathcal{B}$ such that $P \subseteq B_0$.\qed

\medskip
Let $f(x)=\frac{k^2-k}{2}+\frac{x^2-x(2k-3)}{2}.$
It is easy to see that $f(x)$ is a decreasing function in $[2,k-2]$ and $f(k-1)=f(k-2)$. Hence, if $2\leq |P|\leq k-1$ then
\[|I_P|\geq f(|P|)=\frac{k^2-k}{2}+\frac{|P|^2-|P|(2k-3)}{2}.\]

\medskip
\noindent\textbf{Proof of Theorem \ref{NO8}.}
Let $S$ be a dominating set of $D$ with $|S|=\gamma(D)$, $P=\pi(S)$ and $L=S\setminus P$. Since $D$ is a symmetric design, without loss of generality we may assume that $|P|\leq |L|$, for otherwise we can consider the dual design of $D$ whose incidence graph is isomorphic to that of $D$.

Consider the case $k\leq 35$ firstly. Assume on the contrary that
\begin{eqnarray*}
\gamma(D)&\leq&k-2+\sum_{i\geq 1}\left\lfloor \frac{k-4}{2^{2i-1}} \right\rfloor\\
&=&k-2+\left\lfloor \frac{k-4}{2} \right\rfloor+\left\lfloor \frac{k-4}{8} \right\rfloor\\
&\leq& k-2+\left\lfloor \frac{k-4}{2} \right\rfloor+\frac{k-4}{8}=:2s
\end{eqnarray*}
where $s=\frac{1}{2}(k-2+\left\lfloor \frac{k-4}{2} \right\rfloor+\frac{k-4}{8})<k$. Then $|P|\leq \left\lfloor\frac{\gamma(D)}{2}\right\rfloor \leq \frac{\gamma(D)}{2}\leq s<k$.

By Lemma \ref{NO7}, we have $
2s\geq \gamma(D)=|S|\geq |I_P|\geq f(|P|)\geq f(s)$. It follows that $k^2-k+s(s-2k-1)\leq 0$. However, if $k$ is even, then $s=\frac{1}{2}(k-2+\frac{k-4}{2}+\frac{k-4}{8})=\frac{13k-36}{16}$. Thus
$k^2-k+s(s-2k-1)=\frac{9k^2-248k+1872}{256}\geq \frac{164}{256}>0,$
which is a contradiction.
If $k$ is odd, then $s=\frac{1}{2}( k-2+\frac{k-5}{2}+\frac{k-4}{8})=\frac{13k-40}{16}$. Thus
$
k^2-k+s(s-2k-1)=\frac{9k^2-224k+2240}{256}\geq \frac{849}{256}>0,$
which is a contradiction.

Hence $k\geq 36$. Assume on the contrary that
\begin{eqnarray*}
\gamma(D)&\leq&k-2+\sum_{i\geq 1}\left\lfloor \frac{k-4}{2^{2i-1}} \right\rfloor\\
&\leq&k-2+\sum_{i\geq 1}\frac{k-4}{2^{2i-1}} \\
&\leq&k-2+\frac{2}{3}(k-4)\\
&=&\frac{5k-14}{3}=:2t,
\end{eqnarray*}
where $t=\frac{5k-14}{6}<k$. Then $|P|\leq \left\lfloor\frac{\gamma(D)}{2}\right\rfloor \leq \frac{\gamma(D)}{2}\leq t<k$.

By Lemma \ref{NO7}, we have $2t\geq \gamma(D)=|S|\geq |I_P|\geq f(|P|)\geq f(t)$. It follows that
$k^2-k+t(t-2k-1)\leq 0.$ Recall that $t=\frac{5k-14}{6}$ and $k\geq 36$. So we obtain
$$k^2-k+t(t-2k-1)=\frac{k^2-38k+280}{36}\geq\frac{208}{36}>0,$$
a contradiction. This completes the proof of Theorem \ref{NO8}. \qed

\section{Super-neat designs}

\begin{lemma}\label{NO2}
Let $D=(X,\mathcal{B})$ be a $2\text{-}(v,k,\lambda)$ design with $\gamma (D) < \frac{{\left\lceil {\frac{r}{\lambda }}\right\rceil (k-1)+v}}{k}$,
then $D$ is a super-neat design.
In particular, when $\lambda=1$, then $D$ is a super-neat design if $\gamma (D) < \frac{2v-1}{k}$.
\end{lemma}
\noindent\textbf{Proof.}
Assume that $D$ is not a super-neat design. Let $S$ be a dominating but not neat set of $D$ with $|S|=\gamma (D)$. Then
$I_P \subsetneq S$, where $P=\pi (S)$. From the proof of Lemma $\ref{NO1}$, we have $|P|\geq \left\lceil{\frac{r}{\lambda }} \right\rceil$. Now, by Lemma \ref{4.4}, we have $\gamma (D)=|S|\geq \frac{{\left\lceil {\frac{r}{\lambda }}\right\rceil (k-1)+v}}{k}$, which is a contradiction.  \qed

The next theorem proves that finite projective planes are super-neat. This solves Conjecture 6.2 in \cite{Felix}.

\begin{theorem}\label{NO3}
Finite projective planes are super-neat.
\end{theorem}
\noindent\textbf{Proof.} Let $D=(X,\mathcal{B})$ be a finite projective plane of order $q$, that is, a $2\text{-}(q^2+q+1,q+1,1)$ design. By Lemma \ref{5.2}, $\gamma (D)=2q$, and
\[\frac{{2v - 1}}{k}=\frac{2(q^2+q+1)-1}{q+1}=2q+\frac{1}{q+1}>2q=\gamma (D).\]
Hence, by Lemma \ref{NO2}, $D$ is a super-neat design.\qed

\begin{theorem}\label{NO4}
Let $D=(X,\mathcal{B})$ be a finite affine plane of order $q$. Then $\gamma (D)=2q-1$.
\end{theorem}
\noindent\textbf{Proof.}
Here $D$ is a $2$-$(q^2,q,1)$ design. By Theorem \ref{NO6},
\[2q-1\leq \gamma(D)\leq \frac{(v-q^2)(v-1)}{q(q-1)}+2q-1,\]
that is, $2q-1\leq \gamma(D)\leq 2q-1.$
It follows that $\gamma(D)=2q-1$.\qed

\begin{theorem}\label{NO5}
Finite affine planes are super-neat.
\end{theorem}
\noindent\textbf{Proof.} Let $D=(X,\mathcal{B})$ be a finite affine plane of order $q$. By Theorem \ref{NO4}, we have $\gamma (D)=2q-1$. On the other hand, as $q\geq 2$, we have
$$\frac{2v-1}{k}=\frac{2q^2-1}{q}=2q-\frac{1}{q} > 2q-1=\gamma (D).$$
Now the theorem follows from Lemma \ref{NO2}.\qed

\section{Residual designs}

Let $D=(X,\mathcal{B})$ be a symmetric $2\text{-}(v,k,\lambda)$ design and $B_0\in \mathcal{B}$. Then the \textit{residual design} of $D$ is defined by
$Res(X,\mathcal{B},B_0)=(X\setminus B_0,\{B\setminus B_0:B\in \mathcal{B},B\neq B_0\})$.

Let $D_1=(X_1,\mathcal{B}_1)$, $D_2=(X_2,\mathcal{B}_2)$ be two $2\text{-}(v,k,\lambda)$ designs. We say that $D_1$ and $D_2$ are \textit{isomorphic} if there exists a bijection $g:X_1\rightarrow X_2$ such that for any $B_1 \in \mathcal{B}_1$, $B_1^g\in \mathcal{B}_2$ holds. If $D=D_1=D_2$ then $g$ is called an \textit{automorphism} of $D$. All automorphisms of $D$ constitute a group, called the \textit{full automorphism group} of $D$ and denoted by $Aut(D)$. Any subgroup $H$ of $Aut(D)$ is also called an
\textit{automorphism group} of $D$. We say that $H$ is \textit{point-transitive} if $H$ is transitive on $X$, that is for any
$x,y\in X$, there exists $g\in H$ such that $x^g=y$; $H$ is \textit{block-transitive} if $H$ is transitive on $\mathcal{B}$, that is for any $B,C\in
\mathcal{B}$, there exists $g\in H$ such that $B^g=C$.

By the above definitions, the following result is obvious.

\begin{lemma}\label{block}
Let $D=(X,\mathcal{B})$ be a block-transitive symmetric $2\text{-}(v,k,\lambda)$ design.
Then for any $B_1,B_2\in \mathcal{B}$, $
Res(X,\mathcal{B},B_1)\cong Res(X,\mathcal{B},B_2).$
\end{lemma}

\begin{theorem}\label{NO10}
Let $D=(X,\mathcal{B})$ be a symmetric $2\text{-}(v,k,\lambda)$ design and $D_1=Res(X,\mathcal{B},$ $B_0)$, a residual design of $D$, where $B_0\in \mathcal{B}$. Then $\gamma(D_1)\geq\gamma(D)-1$.
\end{theorem}
\noindent\textbf{Proof.}
Let $S$ be a minimum dominating set of $D_1$ so that $|S|=\gamma(D_1)$. Clearly, we have $B_0\notin S$. Let $P=\pi(S)$ and $S'=S\cup\{B_0\}$. Then $S'$ is a dominating set of $D$ because any vertex in $(\mathcal{B}\setminus \{B_0\})\cup (X\setminus B_0)$ is dominated by $S$, and any vertex in $B_0\cup\{B_0\}$ is dominated by $\{B_0\}$. Hence,
$\gamma(D)\leq |S'|=\gamma(D_1)+1$.\qed

\medskip
\noindent\textbf{Proof of Theorem \ref{NORes}.} On the one hand, by Theorem \ref{NO10}, $\gamma(D_1)\geq\gamma(D)-1$.

On the other hand, by Lemma \ref{block}, $\gamma(Res(X,\mathcal{B},B_0))=\gamma (Res(X,\mathcal{B},B_0'))$ for any $B_0,B_0'\in \mathcal{B}$. Let
$S$ be a minimum dominating set so that $|S|=\gamma(D)$. Let $P=\pi(S)$, $L=S\setminus P$ and $B_0\in \hat{L}(P)$. Since any $B\in\mathcal{B}\setminus L$ is dominated by $P$, and any $x\in X\setminus (P\cup B_0)$ is dominated by $L\setminus
\{B_0\}$, then $S\setminus\{B_0\}$ is a dominating
set of $D_1=Res(X,\mathcal{B},B_0)$. Hence $\gamma(D_1)\leq |S\setminus \{B_0\}|=\gamma(D)-1.$

Thus $\gamma(D_1)=\gamma(D)-1$. This completes the proof of Theorem \ref{NORes}. \qed

\subsection*{Acknowledgements}
Thanks to Professor Sanming Zhou at University of Melbourne for corrections and some useful discussion which lead to the improvement of the paper.


\begin{thebibliography}{99}

\bibitem{Hedetniemi}
S. T. Hedetniemi, R. C. Laskar, Bibliography on domination in graphs and some basic definitions of domination parameters, Discrete Math. 86(1-3)(1990) 257-277.

\bibitem{Garey}
M. R. Garey, D. S. Johnson, Computers and Intractability: A Guide to the Theory of NP-Completeness, W. H. Freeman and Company, 1979.

\bibitem{Bondy}
J. A. Bondy, U. S. R. Murty, Graph Theory with Applications,
New York, Macmillan Press, 1976.

\bibitem{Felix}
F. Goldberg, D. Rajendraprasad, R. Mathew, Domination in designs, Arxiv:1405.3436, 2014.

\bibitem{Wan}
Z. X. Wan, Design Theory, Beijing, Higher Education Press, 2009.

\bibitem{Laskar}
R. Laskar and C. Wallis, Chessboard graphs, related designs, and domination parameters, J. Stat.
Plann. Inference 76(1-2)(1999) 285-294.

\bibitem{boll}
B. Bollob\'{a}s, E. Cockayne, Graph-theoretic parameters concerning domination, independence,
and irredundance, J. Graph Theory 3(3)(1979) 241-249.

\end{thebibliography}
\end{document}